\documentclass[12pt]{amsart}
\usepackage{amssymb,latexsym,amsthm}
\usepackage{graphicx}
\usepackage{color}

\newcommand\mgj{M(G_j)}
\newcommand\mgo{M(G_1)}
\newcommand\mgt{M(G_2)}
\newcommand\wgj{\widehat{G_j}}
\newcommand\wgo{\widehat{G_1}}
\newcommand\wgt{\widehat{G_2}}
\newcommand\wgtp{\widehat{G_2}+}
\newcommand\wgtm{\widehat{G_2}-}
\newtheorem{theorem}{Theorem}[section]

\theoremstyle{definition}

\theoremstyle{remark}

\numberwithin{equation}{section}

\begin{document}

\baselineskip=17pt

\title[]
{New criteria for equivalence of locally compact abelian groups}

% author 1
\author{Osamu~Hatori}
\address{Department of Mathematics, Faculty of Science,
Niigata University, Niigata 950-2181 Japan}
\curraddr{}
\email{hatori@math.sc.niigata-u.ac.jp}

\thanks{The  author was partly
supported by the Grants-in-Aid for Scientific
Research, Japan Society for the Promotion of Science.
}



\begin{abstract}
If an open subgroup of the group of the invertible measures on a LCA group is isometric to another one, then the corresponding underlying LCA groups are topologically isomorphic to each other.
\end{abstract}
\maketitle
\noindent
Department of Mathematics,\\ Faculty of Science,
Niigata University,\\ Niigata 950-2181, Japan
\\
email: hatori@math.sc.niigata-u.ac.jp

\newpage
%%%%%%%%%%%%%%%%%%%%%%%%%%%%%%%%%%%%%%%%%%%%%%%%%%%%%%%%%%%%%%%%%%%%%%%%%%%%%%%%%%%%%%%%%%%%%%%%%%%%%%%%%%%%%%%%%%%%%%%%%%%%%%%%%%%%%%%%%%%%%%%%%%%%%%%%%%%%%%%%%%%%%%%%%%%%%%%%%%%%%%%%%%%%%%%%%%%%%%%%%%%%%%%%%%%%%%%%%%%%%%%%%%%%%%%%%%%%%%%%%%%%%%%%%%%%%%%%%%%%%%%%%%%%%%%%%%%%%%%%%%%%%%%%%%%%%%%%%%%%%%%%%%%%%%%%%%%%%%%%%%%%%%%%%%%%%%%%%%%%%%%%%%%%%%%%%%%%%%%%%%%%%%%%%%%%%%%%%%%%%%%%%%%%%%%%%%%%%%%
\section{Introduction}
%%%%%%%%%%%%%%%%%%%%%%%%%%%%%%%%%%%%%%%%%%%%%%%%%%%%%%%%%%%%%%%%%%%%%%%%%%%%%%%%%%%%%%%%%%%%%%%%%%%%%%%%%%%%%%%%%%%%%%%%%%%%%%%%%%%%%%%%%%%%%%%%%%%%%%%%%%%%%%%%%%%%%%%%%%%%%%%%%%%%%%%%%%%%%%%%%%%%%%%%%%%%%%%%%%%%%%%%%%%%%%%%%%%%%%%%%%%%%%%%%%%%%%%%%%%%%%%%%%%%%%%%%%%%%%%%%%%%%%%%%%%%%%%%%%%%%%%%%%%%%%%%%%%%%%%%%%%%%%%%%%%%%%%%%%%%%%%%%%%%%%%%%%%%%%%%%%%%%%%%%%%%%%%%%%%%%%%%%%%%%%%%%%%%%%%%%%%%%%%

A long tradition of inquiry seeks sufficient sets of conditions on morphisms between algebras of measures on locally compact groups in order that the groups are topologically isomorphic to each other. In case the morphism is an algebra isomorphism between measure algebras on locally compact abelian groups (LCA groups), the LCA groups are topologically isomorphic to each other provided that the given isomorphism is also isometric. Kalton and Wood \cite{KW} proved the same conclusion under the additional weaker assumption that the norm of the algebra isomorphism is less than $\sqrt{2}$ and showed that the constant is the best possible. It had already been known that the measure algebras are algebraically isomorphic to each other if and only if there exists a piecewise affine homeomorphism between the dual groups of the underlying LCA groups (\cite[4.6.4,\,\,4.7.7]{R}), and in fact there are several non-isomorphic LCA groups between the dual groups of which piecewise affine homeomorphisms exist (cf. \cite[Theorem in 4.7.7]{R}). Furthermore measure algebras on any discrete groups with the same cardinality are isometrically isomorphic as Banach spaces to each other. These circumstances reveal that the algebraic structure nor the structure as a Banach space of the measure algebra {\it alone} do not encode completely in the group structure of the underlying LCA groups.

In this paper we propose criteria for the equivalence of LCA groups from another point of view; the metrical structure of an open subgroup of the group of the invertible measures in the measure algebra. We show as the main result in this paper that the coincidence of the metric structures of those open subgroups guarantees the coincidence of the structure of the underlying LCA groups.

Throughout the paper $G_j$ is a locally compact abelian group for $j=1,2$ with the group operation denoted by the symbol $+$ and the unit element $0$. The measure algebra on $G_j$ is denoted by $M(G_j)$ and the ideal of all absolutely continuous measures with respect to the Haar measure is $L^1(G_j)$. Note that $L^1(G_j)$ is identified with the algebra of all complex-valued integrable functions with respect to the Haar measure on $G_j$. The continuous character on $G_j$ is denoted by $(\cdot,\gamma )$ and the dual group of $G_j$ is denoted by $\widehat{G_j}$. For a measure $\mu$ in $\mgj$, the Fourier-Stieltjes transform of $\mu$ is denoted by $\hat{\mu}$ and $\widehat{\mgj}=\{\hat{\mu}:\mu \in \mgj \}$. It is well-known that $\mgj$ is a unital semisimple commutative Banach algebra. The Gelfand transform of $\mu \in \mgj$ is denoted by $\Gamma(\mu)$. The dual group $\widehat{G_j}$ is embedded naturally in $\Phi_{\mgj}$ and $\Gamma(\mu)=\hat{\mu}$ for every $\mu\in \mgj$ (cf. \cite[p.309]{LN}). For an $x\in G_j$, $\delta_x$ denotes the Dirac measure at $x$.

%%%%%%%%%%%%%%%%%%%%%%%%%%%%%%%%%%%%%%%%%%%%%%%%%%%%%%%%%%%%%%%%%%%%%%%%%%%%%%%%%%%%%%%%%%%%%%%%%%%%%%%%%%%%%%%%%%%%%%%%%%%%%%%%%%%%%%%%%%%%%%%%%%%%%%%%%%%%%%%%%%%%%%%%%%%%%%%%%%%%%%%%%%%%%%%%%%%%%%%%%%%%%%%%%%%%%%%%%%%%%%%%%%%%%%%%%%%%%%%%%%%%%%%%%%%%%%%%%%%%%%%%%%%%%%%%%%%%%%%%%%%%%%%%%%%%%%%%%%%%%%%%%%%%%%%%%%%%%%%%%%%%%%%%%%%%%%%%%%%%%%%%%%%%%%%%%%%%%%%%%%%%%%%%%%%%%%%%%%%%%%%%%%%%%%%%%%%%%%%
\section{A form of real-algebra isomorphisms}
%%%%%%%%%%%%%%%%%%%%%%%%%%%%%%%%%%%%%%%%%%%%%%%%%%%%%%%%%%%%%%%%%%%%%%%%%%%%%%%%%%%%%%%%%%%%%%%%%%%%%%%%%%%%%%%%%%%%%%%%%%%%%%%%%%%%%%%%%%%%%%%%%%%%%%%%%%%%%%%%%%%%%%%%%%%%%%%%%%%%%%%%%%%%%%%%%%%%%%%%%%%%%%%%%%%%%%%%%%%%%%%%%%%%%%%%%%%%%%%%%%%%%%%%%%%%%%%%%%%%%%%%%%%%%%%%%%%%%%%%%%%%%%%%%%%%%%%%%%%%%%%%%%%%%%%%%%%%%%%%%%%%%%%%%%%%%%%%%%%%%%%%%%%%%%%%%%%%%%%%%%%%%%%%%%%%%%%%%%%%%%%%%%%%%%%%%%%%%%%

While the representation through the Gelfand transform of {\it complex} isomorphism between unital semisimple commutative complex Banach algebras is quite popular, the author first remark that he was unable to find the description of the general form of surjective real-algebra isomorphisms between those algebras. We begin with a statement on the form of the Gelfand transform. For the completeness we also give a proof of the statement.

%%%%%%%%%%%%%%%%%%%%%%%%%%%%%%%%%%%%%%%%%%%%%%%%%%%%%%%%%%%%%%%%%%%%%%%%%%%%%%%%%%%%%%%%%%%%%%%%%%%%%%%%%%%%%%%%%%%%%%%%%%%%%%%%%%%%%%%%%%%%%%%%%%%%%%%%%%%%%%%
\begin{theorem}\label{gf}
Let $A$ and $B$ be unital semisimple commutative (complex) Banach algebras with the unit $e_A$ and $e_B$ respectively. Suppose that $S$ is a surjective real-algebra isomorphism from $A$ onto $B$. Then there exists a homeomorphism $\Psi$ from $\Phi_B$ onto $\Phi_A$ such that the Gelfand transform $\Gamma (S(a))$ for $a\in A$ has the form
\begin{equation}
\Gamma (S(a))(\phi)=
\begin{cases}
\Gamma (a)\circ \Psi (\phi), & \phi \in \Phi_{B+},
\\
\\
\overline{
\Gamma (a)\circ \Psi (\phi)
}, &\phi \in \Phi_B\setminus \Phi_{B+},
\end{cases}
\end{equation}
where $\Phi_{B+}=\{\phi \in \Phi_B:\Gamma(S(ie_A)(\phi)=i\}$.
\end{theorem}
%%%%%%%%%%%%%%%%%%%%%%%%%%%%%%%%%%%%%%%%%%%%%%%%%%%%%%%%%%%%%%%%%%%%%%%%%%%%%%%%%%%%%%%%%%%%%%%%%%%%%%%%%%%%%%%%%%%%

%%%%%%%%%%%%%%%%%%%%%%%%%%%%%%%%%%%%%%%%%%%%%%%%%%%%%%%%%%%%%%%%%%%%%%%%%%%%%%
\begin{proof}
To begin the proof, first observe that $S(e_A)=e_B$ is rather trivial since $S$ is surjective and $\Gamma$ is injective. It is easy to check that $\Gamma (S(ie_A))$ takes the values $i$ or $-i$ for $(S(ie_A))^2=-e_B$. Therefore the sets $\Phi_{B+}$ and $\Phi_{B-}=\{\phi\in \Phi_B:\Gamma (S(ie_A))(\phi)=-i\}$ are possibly empty disjoint closed and open subsets of $\Phi_B$, as $\Gamma (S(ie_A))$ is continuous, and $\Phi_{B+}\cup \Phi_{B-}=\Phi_B$. Note that $\Phi_B\setminus \Phi_{B+}=\Phi_{B-}$.

To define the function $\Psi:\Phi_B\to \Phi_A$, consider $\Psi(\phi)$ in two cases: $\phi\in \Phi_{B+}$; $\phi\in \Phi_{B-}$. First, let $\phi \in \Phi_{B+}$. The functional $\phi\circ S$ on $A$ is  multiplicative and real-linear, and also non-trivial for $\phi(S(e_A))=1$. It requires only easy computation to verify that $\phi\circ S(ia)=\phi\circ S(ie_A)\phi \circ S(a)=i\phi \circ S(a)$ holds for every $a\in A$, ensuring that $\phi \circ S$ is complex-linear as is already real-linear. Next considering the second case where $\phi\in \Phi_{B-}$, it can be checked in a way similar to the first case that $\overline{\phi\circ S}$ gives a non-trivial multiplicative complex linear functional on $A$. Therefore
\begin{equation}\label{circ}
\Psi (\phi)=
\begin{cases}\phi\circ S, & \phi \in \Phi_{B+} \\
\\
\overline{\phi\circ S}, & \phi \in \Phi_{B-}
\end{cases}
\end{equation}
is a well-defined map from $\Phi_B$ into $\Phi_A$. Rewriting the equation (\ref{circ}) we obtain for each $a\in A$
\begin{equation}\label{circcirc}
\phi(S(a))=
\begin{cases}
\Psi(\phi)(a), &\text{if $\phi \in \Phi_{B+}$}
\\
\\
\overline{\Psi (\phi)(a)}, &\text{if $\phi \in \Phi_{B-}$}
\end{cases}
\end{equation}
and due to the definition of the Gelfand transform to get
\begin{equation}\label{gamma}
\Gamma (S(a))(\phi )= 
\begin{cases}
\Gamma (a)\circ\Psi (\phi), &\text{if $\phi \in \Phi_{B+}$}
\\
\\
\overline{\Gamma(a)\circ \Psi (\phi)}, &\text{if $\phi \in \Phi_{B-}$},
\end{cases}
\end{equation}
which is the desired form. 

We assert $\Psi$ is the desired homeomorphism from $\Phi_B$ onto $\Phi_A$. First observe that $\Psi$ is surjective. Let $\varphi \in \Phi_A$. We will find out $\tilde \varphi \in \Phi_B$ with $\Psi(\tilde\varphi)=\varphi$. Applying the \v Silov idempotent theorem there is $v\in B$ such that $\Gamma (v)=1$ on $\Phi_{B+}$ and $\Gamma (v)=0$ on $\Phi_{B-}$. By a simple calculation $\Gamma (S(ie_A)v)=\Gamma (S(ie_A))\Gamma (v)=i\Gamma (v)$ hold as $\Gamma (v)=1$ on $\Phi_{B+}$ and $\Gamma (v)$ vanishes on $\Phi_{B-}$, hence $S(ie_A)v=iv$ holds for $\Gamma$ is injective. We obtain
\begin{equation}\label{star}
iS^{-1}(v)=S^{-1}(S(ie_A))S^{-1}(v)=S^{-1}(S(ie_A)v)=S^{-1}(iv),
\end{equation}
as $S^{-1}$ is multiplicative. Since $v^2=v$ holds and also $S^{-1}$ is multiplicative, $\varphi (S^{-1}(v))=1$ or $0$. 

Suppose that $\varphi (S^{-1}(v))=1$. Put $\tilde\varphi:B\to {\mathbb C}$ as $\tilde\varphi (b)=\varphi(S^{-1}(b))$ for each $b\in B$. Then $\tilde \varphi$ is real-linear, multiplicative and $\tilde\varphi(e_B)=1$. Applying the equation (\ref{star})
\begin{multline*}
\tilde\varphi(ib)=\varphi(S^{-1}(ib))\varphi(S^{-1}(v))=\varphi(S^{-1}(ivb)) 
=\varphi(S^{-1}(iv))\varphi(S^{-1}(b))\\
=\varphi(iS^{-1}(v))\tilde\varphi (b)=i\varphi (S^{-1}(v))\tilde\varphi (b)=i\tilde\varphi (b)
\end{multline*}
holds for every $b\in B$. This implies that $\tilde\varphi\in \Phi_B$. It requires only elementary calculation to check that $\tilde\varphi (S(ie_A))=i$, hence $\tilde\varphi \in \Phi_{B+}$. For every $a\in A$
\[
\Psi (\tilde\varphi )(a)=\tilde\varphi\circ S(a)=\varphi (a)
\]
holds, hence $\Psi (\tilde\varphi )=\varphi$. 

Suppose next that $\varphi(S^{-1}(v))=0$. In this case $\varphi (S^{-1}(e_A-v))=1$ holds. Put $\tilde\varphi:B\to {\mathbb C}$ as $\tilde\varphi(b)=\overline{\varphi(S^{-1})b)}$ for $b\in B$ in this case. Then $\tilde \varphi$ is multiplicative, real-linear and $\tilde\varphi(e_B)=1$. Since $\Gamma (e_b-v)=1$ on $\Phi_{B-}$ and $\Gamma (e_B-v)$ vanishes on $\Phi_{B+}$, $\Gamma (S(ie_A)(e_B-v))=-i\Gamma (e_B-v)$ holds, so that $S(ie_A)(e_B-v)=-i(e_B-v)$ holds. Hence
\begin{multline*}
iS^{-1}(e_B-v)=S^{-1}(S(ie_A)(e_B-v))\\
=S^{-1}(-i(e_B-v))=-S^{-1}(i(e_B-v))
\end{multline*}
hold, whence 
\begin{multline*}
\tilde\varphi(ib)=\overline{\varphi (S^{-1}(ib))}
=\overline{\varphi(S^{-1}(e_B-v))\varphi(S^{-1}(ib))}\\
=\overline{\varphi(S^{-1}(i(e_B-v))\varphi(S^{-1}(b))} 
=\overline{\varphi (-iS^{-1}(e_B-v))}\tilde\varphi (b)\\
=i\overline{\varphi(S^{-1}(e_B-v))}\tilde\varphi (b)=i\tilde\varphi (b)
\end{multline*}
hold for every $b\in B$. It turns that $\tilde\varphi$ is complex-linear and $\tilde\varphi\in \Phi_B$. We also have
\[
\tilde\varphi (S(ie_A))=\overline{\varphi(S^{-1}(S(ie_A)))}=\overline{\varphi (ie_A)}=-i,
\]
hence $\tilde\varphi \in \Phi_{B-}$. Furthermore
\[
\Psi(\tilde\varphi)(a)=\overline{\tilde\varphi\circ S(a)}=\varphi (a)
\]
holds for every $a\in A$, hence $\Psi (\tilde\varphi)=\varphi$. With the observation for the case where $\varphi(S^{-1}(v))=1$ we conclude that $\Psi$ is surjective.

We claim that $\Psi$ is injective. Let $\phi_1$ and $\phi_2$ be different points in $\Phi_B$. Then there exists $b\in B$ with $\phi_1(b)=0$ and $\phi_2(b)=1$. It requires only a simple calculation due to the definition of $\Psi$ that $\Psi(\phi_1)(S^{-1}(b))\ne\Psi (\phi_2)(S^{-1}(b))$, whence $\Psi (\phi_1)\ne\Psi(\phi_2)$. 

A proof of the continuity of $\Psi$ is as routine as follows. Suppose that $\phi\in \Phi_B$ and a net $\{\phi_{\alpha}\}$in $\Phi_B$ converges to $\phi$. The two sets $\Phi_{B+}$ and $\Phi_{B-}$ are closed and open, so there is no serious loss of generality we may assume that $\{\phi_{\alpha}\}\subset \Phi_{B+}$ if $\phi \in \Phi_{B+}$ and $\{\phi_{\alpha}\}\subset \Phi_{B-}$ if $\phi\in \Phi_{B-}$. Then 
\[
\Psi (\phi_{\alpha})(a)=\phi_{\alpha}(S(a))\to \phi(S(a))=\Psi(\phi)(a)
\]
if $\phi\in \Phi_{B+}$ and 
\[
\Psi(\phi_{\alpha})(a)=\overline{\phi_{\alpha}(S(a))}\to \overline{\phi_{\alpha}(S(a))}=\Psi(\phi)(a)
\]
if $\phi\in \Phi_{B-}$ hold for every $a\in A$. We conclude that the net $\{\Psi(\phi_{\alpha})\}$ converges to $\Psi (\phi)$ since the Gelfand topology on $\Phi_B$ is the weak topology induced by $\Gamma(B)$. Hence $\Psi$ is continuous. We assert that the injective continuous surjection $\Psi$ is a homeomorphism as $\Phi_B$ is compact and $\Phi_A$ is a Hausdorff space.
\end{proof}
%%%%%%%%%%%%%%%%%%%%%%%%%%%%%%%%%%%%

%%%%%%%%%%%%%%%%%%%%%%%%%%%%%%%%%%%%%%%%%%%%%%%%%%%%%%%%%%%%%%%%%%%%%%%%%%%%%%%%%%%%%%%%%%%%%%%%%%%%%%%%%%%%%%%%%%%%%%%%%%%%%%%%%%%%%%%%%%%%%%%%%%%%%%%%%%%%%%%%%%%%%%%%%%%%%%%%%%%%%%%%%%%%%%%%%%%%%%%%%%%%%%%%%%%%%%%%%%%%%%%%%%%%%%%%%%%%%%%%%%%%%%%%%%%%%%%%%%%%%%%%%%%%%%%%%%%%%%%%%%%%%%%%%%%%%%%%%%%%%%%%%%%%%%%%%%%%%%%%%%%%%%%%%%%%%%%%%%%%%%%%%%%%%%%%%%%%%%%%%%%%%%%%%%%%%%%%%%%%%%%%%%%%%%%%%%%%%%%
\section{The main result}
%%%%%%%%%%%%%%%%%%%%%%%%%%%%%%%%%%%%%%%%%%%%%%%%%%%%%%%%%%%%%%%%%%%%%%%%%%%%%%%%%%%%%%%%%%%%%%%%%%%%%%%%%%%%%%%%%%%%%%%%%%%%%%%%%%%%%%%%%%%%%%%%%%%%%%%%%%%%%%%%%%%%%%%%%%%%%%%%%%%%%%%%%%%%%%%%%%%%%%%%%%%%%%%%%%%%%%%%%%%%%%%%%%%%%%%%%%%%%%%%%%%%%%%%%%%%%%%%%%%%%%%%%%%%%%%%%%%%%%%%%%%%%%%%%%%%%%%%%%%%%%%%%%%%%%%%%%%%%%%%%%%%%%%%%%%%%%%%%%%%%%%%%%%%%%%%%%%%%%%%%%%%%%%%%%%%%%%%%%%%%%%%%%%%%%%%%%%%%%%
For a measure algebra $M(G)$ on a LCA group $G$, $M(G)^{-1}$ denotes the group of all invertible elements in the algebra $M(G)$. If $G_1$ and $G_2$ are discrete LCA groups with the same cardinal number and $\pi$ is {\it any} bijection from $G_2$ onto $G_1$. In this case $M(G_j)$ coincides with $L^1(G_j)$ and $f \to f\circ\pi$ defines a isometrical linear isomorphism from $L^1(G_1)$ onto $L^1(G_2)$. In general this linear isomorphism need not be an {\it algebra} isomorphism. The structure as a Banach space of the measure algebra does not encode in the group structure of the underlying LCA group, although the structure as a Banach algebras does. Since the invertible elements are clearly related to the multiplication, and the metric structure of an open sets encodes the linear structure of the algebra (cf. the local Mazur-Ulam theorem), there arises the natural question of whether and how the metric structure of an open subgroup of $M(G)^{-1}$ encodes the group structure of $G$. The following is the main result in this paper.
%%%%%%%%%%%%%%%%%%%%%%%%%%%%%%%%%%%%%%%%%%%%%%%%%%%%%%%%%%%%%%%%%%%%%%%%%%%%%%%%%%%%%%%%%%%%%%%%%%%%%%%%%%%%%%%%%%%%%%%%%%%%%%%%%%%%%%%%%%%%%%%%%%%%%%%%%%%%%%%
\begin{theorem}\label{img}
Let ${\mathcal A}_j$ be an open subgroup of $M(G_j)^{-1}$ for $j=1,2$. Suppose that $T$ is a surjective isometry from ${\mathcal A}_1$ onto ${\mathcal A}_2$. Then $T$ is extended to a surjective either complex or conjugate linear isometry, say $T_U$, from $M(G_1)$ onto $M(G_2)$ and $T(\delta_0)^{-1}T_U$ is an either complex or conjugate algebra isomorphism. Furthermore $G_1$ is topologically isomorphic to $G_2$ as a topological group.
\end{theorem}
%%%%%%%%%%%%%%%%%%%%%%%%%%%%%%%%%%%%%%%%%%%%%%%%%%%%%%%%%%%%%%%%%%%%%%%%%%%%%%%%%%%%%%%%%%%%%%%%%%%%%%%%%%%%%%%%

\begin{proof}
Since $M(G_j)$ is a unital commutative Banach algebra which is also semisimple \cite[p.309]{LN}, $T$ is extended to a real linear isometry $T_U$ from $M(G_1)$ onto $M(G_2)$ by \cite[Theorem 3.1]{hatori}. Denoting $S=T(\delta_0)^{-1}T_U$ it is a surjective isometric real algebra isomorphism by \cite[Theorem 3.3]{hatori}. Due to Theorem \ref{gf} there exists a homeomorphism $\Psi$ from $\Phi_{M(G_2)}$ onto $\Phi_{M(G_1)}$ such that the Gelfand transform $\Gamma (S(\mu))$ for $\mu \in M(G_1)$ has the form 
\begin{equation}\label{1}
\Gamma (S(\mu ))(\phi)= 
\begin{cases}
\Gamma (\mu)\circ \Psi (\phi), & \text{if $\phi \in K$}
\\
\\
\overline{
\Gamma (\mu )\circ \Psi (\phi )}, & \text{if $\phi \in \Phi_{M(G_2)}\setminus K$},
\end{cases}
\end{equation}
where $K=\{\phi\in \Phi_{M(G_2)}:\Gamma (S(i\delta_0))(\phi)=i\}$. Note that $\Gamma (S(i\delta_0))$ takes the values $i$ or $-i$ since $(\Gamma (S(i\delta_0)))^2=-\Gamma (\delta_0)=-1$ on $\Phi_{M(G_2)}$ and hence $\Phi_{M(G_2)}\setminus K=\{\phi\in \Phi_{M(G_2)}:\Gamma (S(i\delta_0))(\phi)=-i\}$. It is well-known that $\wgj \subset \Phi_{\mgj}$ and $\Gamma (\nu)=\hat{\nu}$ on $\wgj$ for every $\nu \in \mgj$ (cf. \cite[p.309]{LN}). Letting $\wgtp=K\cap\wgt$ and $\wgtm=(\Phi_{M(G_2)}\setminus K)\cap \wgt$ we have
\begin{equation}\label{2}
\Gamma (\mu)\circ \Psi (\phi)=
\begin{cases}
\widehat{S(\mu)}(\phi), &\text{if $\phi \in \wgtp$}
\\
\\
\overline{\widehat{S(\mu)}(\phi)},&\text{if $\phi\in \wgtm$}
\end{cases}
\end{equation}
for every $\mu\in \mgo$ and $\phi\in \wgt$ by rewriting the equation (\ref{1}).

We observe that $\wgtp$ or $\wgtm$ is empty. Assume contrary that $\wgtp\ne\emptyset$ and $\wgtm\ne \emptyset$ we will arrive at a contradiction. Put
\[
J=\frac12 (-iS(i\delta_0)+\delta_0),
\]
where the firstly appeared $\delta_0$ in the right-hand side of the equation denotes the Dirac measure on $G_1$ at the identity elements $0$ of $G_1$ while the second $\delta_0$ is the Dirac measure on $G_2$ at the identity element of $0$ of $G_2$. 
It is easy to check that $J$ is the idempotent measure in $M(G_2)$  in the sense that $J^2=J$ and that $\|J\|\le 1$ since $S$ is an isometry. By the theory of commutative Banach algebras the inequality $\|\hat J\|_{\infty} \le \|J\|$ holds for the supremum norm $\|\cdot\|_{\infty}$ on $\Phi_B$, and $\hat J=1$ on $\wgtp$ ensures that $\|J\|=1$. 

By a simple computation due to the definition of $J$ the set $\{\phi\in\wgt:\hat J(\phi)=1\}$ coincides with $\wgtp$, and is an open coset in $\wgt$ by \cite[3.2.4]{R}. In the same way applying $\delta_0-J$ instead of $J$ we see that $\wgtm$is also an open coset in $\wgt$. These observations imply that there are two possible cases: i) $0\in \wgtp$, that is, $\wgtp$ is an open subgroup ; ii) $0\in \wgtm$, that is, $\wgtm$ is an open subgroup. Suppose $0\in \wgtp$. There is an open subgroup $H$ of $\wgt$ and $\phi\in \wgtm$ with $\wgt-=\phi +H$. Since $(\phi +\wgtp)\cap \wgtp=\emptyset$, $\phi+\wgtp\subset \wgtm=\phi+H$ hold, hence $\wgtp \subset H$ is obtained. On the other hand $\wgtm\cap H=\emptyset$ for $0\not\in \wgtm$, so that $H\subset \wgtp$ and $H=\wgtp$ holds. We conclude that $\wgtm =\phi + \wgtp$ for a $\phi\in \wgtm$ if $0\in \wgtp$. In the same way, if $0\in \wgtm$, then $\wgtp=\phi'+\wgtm$ for a $\phi'\in \wgtp$. 

We claim that there exists and $x_1\in G_2$ such that 
\begin{equation*}
\widehat{\delta_{x_1}}(\phi )= 
\begin{cases}
1, & \text{if $\phi \in \wgtp$}
\\
-1,& \text{if $\phi\in \wgtm$}.
\end{cases}
\end{equation*}
We show only for the case where $0\in \wgtp$ but a proof for the case where $0\in \wgtm$ is similar. We have already learnt that $\wgtm$ is open and $\wgtm=\phi+\wgtp$. Therefore $\wgtp$ is a closed subgroup of $\wgt$ and 
\[
\wgt /\wgtp \cong {\mathbb Z}_2,
\]
where ${\mathbb Z}_2={\mathbb Z}/2{\mathbb Z}$ for the discrete group of the integers ${\mathbb Z}$. Then 
\begin{equation}\label{ks}
\left(\wgt/\wgtp\right)^{\widehat{\,\,\,}} \cong \widehat{{\mathbb Z}_2}={\mathbb Z}_2.
\end{equation}
Applying the celebrated duality theorem of Pontryagin the annihilator $\Lambda$ of $\wgtp$ is $\{x\in G_2 : \text{$(x,\phi)=1$ for $\forall \phi \in \wgtp$}\}$ and it is isomorphic to $(\wgt/\wgtp)^{\widehat{\,\,\,}}$ by \cite[Theorem 2.1.2]{R}, and is isomorphic to ${\mathbb Z}_2$ by (\ref{ks}). Then $\Lambda$ is a two-point set. Let $x_1$ be the unique non-zero element of $\Lambda$.  Since 
\[
\widehat{\delta_{x_1}}(\varphi )=\int (-x,\varphi)d\delta_{x_1}(x)=\overline{(x_1,\varphi)}, \quad \varphi\in \wgt,
\]
$\widehat{\delta_{x_1}}(\varphi)=1$ if and only if $(x_1,\varphi )=1$ for $\varphi\in \wgt$. Since $(0,\varphi)=1$ for any $\varphi \in \wgt$, $\phi \in \wgt$ is in the annihilator of the annilator $\Lambda$ if and only if $\widehat{\delta_{x_1}}(\phi)=1$. As the annihilator of $\Lambda$ is $\wgtp$ itself by \cite[Lemma 2.1.3]{R}, we conclude that $\widehat{\delta_{x_1}}(\phi)=1$ is and only if $\phi \in \wgtp$. 
On the other hand $\widehat{\delta_{x_1}}$ takes the value $1$ and $-1$ for $(\delta_{x_1})^2=\delta_{x_1+x_1}=\delta_0$. It follows that 
\begin{equation}\label{01}
\widehat{\delta_{x_1}}(\phi )= 
\begin{cases}
1, & \text{if $\phi\in \wgtp$}
\\
-1,&\text{if $\phi \in \wgtm$}.
\end{cases}
\end{equation}
Put $\nu_{01} =\delta_0+i\delta_{x_1}\in \mgt$ and $\mu_{01}=S^{-1}(\nu_{01})$. Then $\|\mu_{01}\|=\|S^{-1}(\nu_{01})\|=\|\nu_{01}\|=2$ for $S$ is a real-linear isometry and the norm of the measure is defined by the total variation. Substituting this $\mu_{01}$ in the equation (\ref{2}) we have
\begin{equation*}
\Gamma (\mu_{01})\circ\Psi (\phi)=
\begin{cases}
\widehat{\nu_{01}}(\phi),& \text{if $\phi \in \wgtp$}
\\
\overline{\widehat{\nu_{01}}(\phi)},&\text{if $\phi \in \wgtm$}
\end{cases}
=1+i
\end{equation*}
by the equations $\widehat{\delta_0}=1$ on $\wgt$ and (\ref{01}).
If $\Psi (\wgt)\supset \wgo$ will be claimed, then $\hat{\mu_{01}}=1+i$ on $\wgo$ will follow and it will imply $\mu_{01}=(1+i)\delta_0$ by the uniqueness theorem for the Fourier-Stieltjes transform, which will be a contradiction since $\|\mu_{01}\|=\|(1+i)\delta_0\|=\sqrt{2}$. To prove $\Psi (\wgt)\supset \wgo$ we now need the inclusion $S(L^1(G_1))\subset L^1(G_2)$. For the proof of the inclusion $S(L^1(G_1))\subset L^1(G_2)$  define a deformed map $\check{S}$ of $S$ and show that $\check{S}(L^1(G_1))=L^1(G_2)$ first. 
For a measure $\mu\in \mgj$, let the measure $\tilde{\mu}$ in $\mgj$ be defined as $\tilde{\mu}(E)=\overline{\mu(-E)}$ for every measurable subset $E$ of $G_j$ (cf. \cite[p.105]{R}). Note that $\widehat{\tilde{\mu}}=\overline{\hat{\mu}}$ on $\wgj$. Recall that $J=\frac12(-iS(i\delta_0)+\delta_0)$ is an idempotent. Since $\hat{J}$ is a real-valued function the equations $\widehat{\tilde{J}}=\overline{\hat{J}}=\hat{J}$ holds, hence $\tilde{J}=J$ by the uniqueness theorem for Fourier-Stieltjes transforms. 
Let the map $\check{S}:\mgo \to \mgt$ be defined as $\check{S}(\mu)=JS(\mu)+\widetilde{((\delta_0-J)(S(\mu))}$ for each $\mu\in \mgo$. Recall that $\hat J=1$ on $\wgtp$ and $\hat J=0$ on $\wgtm$. Together with the observation that $\widehat{\widetilde{((\delta_0-J)S(\mu))}}$ is the complex conjugate of $\widehat{(\delta_0-J)S(\mu)}$ on $\wgt$ and the equation (\ref{1}) we have by a simple computation that 
\[
\widehat{\check{S}(\mu)}(\phi)=\Gamma (\mu)\circ \Psi (\phi),\quad \mu \in \mgo,\,\,\phi\in \wgt
\]
holds. By the uniqueness theorem for the Fourier-Stieltjes transform $\check{S}$ is a complex algebra homomorphism. Furthermore $\check{S}$ is injective and surjective. A proof is as follows. Let $\mu,\mu'\in \mgo$ and $\check{S}(\mu)=\check{S}(\mu')$. If $\phi\in \wgtp$, then 
\[
\widehat{\check{S}(\mu)}(\phi)=\widehat{JS(\mu)}(\phi)=\widehat{S(\mu)}(\phi),
\]
\[
\widehat{\check{S}(\mu')}(\phi)=\widehat{JS(\mu')}(\phi)=\widehat{S(\mu')}(\phi),
\]
for $\hat J=1$ and $(\delta_0-J)^{\widehat{\,\,\,}}=0$ on $\wgtp$. If $\phi \in \wgtm$, then 
\[
\widehat{\check{S}(\mu)}(\phi)=\overline{\widehat{((\delta_0-J)S(\mu))}(\phi)}=\overline{\widehat{S(\mu)}(\phi)},
\]
\[
\widehat{\check{S}(\mu')}(\phi)=\overline{\widehat{((\delta_0-J)S(\mu'))}(\phi)}=\overline{\widehat{S(\mu')}(\phi)},
\]
Therefore $\widehat{S(\mu)}=\widehat{S(\mu')}$ on $\wgt$, hence $S(\mu)=S(\mu')$ and $\mu=\mu'$ for $S$ is injective. This implies that $\check{S}$ is injective. To prove the surjectivity of $\check{S}$ let $\nu\in \mgt$. Put $\check{\nu}=J\nu+\widetilde{((\delta_0-J)\nu)}$ and $\check{\mu}=S^{-1}(\check{\nu})$. Then 
\begin{multline*}
\check{S}(\check{\mu})=JS(\check{\mu})+\widetilde{((\delta_0-J)S(\check{\mu}))}=J\check{\nu}+\widetilde{((\delta_0-J)(\check{\nu}))} \\
=J(J\nu + \widetilde{((\delta_0-J)\nu)})+(\delta_0-J)\widetilde{(J\nu+\widetilde{((\delta_0-J)\nu)})}=\nu.
\end{multline*}
hold and the surjectivity of $\check{S}$ is observed. The above observations imply that $\check{S}$ is a complex algebra isomorphism from $\mgo$ onto $\mgt$. Applying Theorem 4.6.4 in \cite{R} $\check{S}(L^1(G_1))=L^1(G_2)$ is derived. 
It is easy to check by a simple calculation that $S=J\check{S}+\widetilde{((\delta_0-J)\check{S})}$ holds. Since $L^1(G_2)$ is an ideal of $M(G_2)$ and closed under the operation $\tilde{\cdot}$ we have the inclusion $S(L^1(G_1))\subset L^1(G_2)$. 

Now we prove that $\Psi(\wgt)\supset \wgo$ by applying the inclusion above. Let $\varphi\in \wgo$. Since $\Psi(\Phi_{\mgt})=\Phi_{\mgo}$ and $\wgo\subset \Phi_{\mgo}$, there is $m\in \Phi_{\mgt}$ with $\Psi (m)=\varphi$. Choose an $f\in L^1(G_1)$ such that $\hat{f}(\varphi)\ne0$. Such an $f\in L^1(G_1)$ exists for $\varphi\in \wgo$. Then $S(f)$ is in $L^1(G_2)$ and the equation (\ref{1})
\begin{multline*}
\Gamma (S(f))(m)=
\begin{cases}
\Gamma (f)\circ \Psi (m), & \text{if $m\in K$}
\\
\overline{\Gamma (f)\circ \Psi (m)}, &\text{if $m\in \Phi_{\mgt}\setminus K$}
\end{cases}
\\
=
\begin{cases}
\hat{f}(\varphi),&\text{if $m\in K$} 
\\
\overline{\hat{f}(\phi)}, & \text{if $m\in \Phi_{\mgt} \setminus K$}.
\end{cases}
\end{multline*}
hold. The above second equality follows by the equation $\Gamma (f)=\hat{f}$ on $\wgo$. Hence $\Gamma (S(f))(m)\ne 0$. Since $\Gamma (g)=0$ on $\Phi_{\mgt}\setminus \wgt$ for every $g\in L^1(G_2)$ (cf. \cite[p.365]{LN}), we obtain that $m\in \wgt$ for $S(f)$ is in $L^1(G_2)$. Therefore we conclude that $\Psi (\wgt)\supset \wgo$ hold as is desired. As is already pointed out before we soon arrive at a contradiction proving that one of $\wgtp$ or $\wgtm$ is empty.

Suppose that $\wgtm$ is empty. Considering $\phi \in \wgt$ in the equation (\ref{1}) $\widehat{S(\mu)}(\phi)=\Gamma (\mu)\circ\Psi (\phi)$ follows. Hence $S$ preserves the complex scalars by the uniqueness theorem for the Fourier-Stieltjes transform. We conclude in this case that $T(\delta_0)^{-1}T_U$ is an isometrical complex algebra isomorphism and $G_1$ is topologically isomorphic to $G_2$ by Theorems 4.6.4 and 4.7.1 in \cite{R}.

Suppose that $\wgtp$ is empty. Then as in a similar way as above $T(\delta_0)^{-1}T_U$ is an isometrical conjugate algebra isomorphism. Then $\widetilde{(T(\delta_0)^{-1}T_U)}$ gives an isometrical complex algebra isomorphism form $\mgo$ onto $\mgt$, hence $G_1$ is topologically isomorphic to $G_2$ as before.
\end{proof}

Acknowledgement. The author would like to express his hearty gratitude toward Professor Jun-ichi Tanaka for suggesting the applications of structure theorems for isometries between groups of invertible elements in unital Banach algebras to a problem on equivalence for locally compact abelian groups.
%%%%%%%%%%%%%%%%%%%%%%%%%%%%%%%%%%%%%%%%%%%%%%%%%%%%%%%%%%%%%%%%%%%%%%%%%%%%%%%%%%%%%%%%%%%%%%%%%%%%%%%%%%%%%%%%%%%%%%%%%%%%%%%%%%%%%%%%%%%%%%%%%%%%%%%%%%%%%%%%%%%%%%%%%%%%%%%%%%%%%%%%%%%%%%%%%%%%%%%%%%%%%%%%%%%%%%%%%%%%%%%%%%%%%%%%%%%%%%%%%%%%%%%%%%%%%%%%%%%%%%%%%%%%%%%%%%%%%%%%%%%%%%%%%%%%%%%%%%%%%%%%%%%%%%%%%%%%%%%%%%%%%%%%%%%%%%%%%%%%%%%%%%%%%%%%%%%%%%%%%%%%%%%%%%%%%%%%%%%%%%%%%%%%%%%%%%%%%%%%%%%%%%%%%%%%%%%%%%%%%%%%%%%%%%%%%%%%%%%%%%%%%%%%%%%%%%%%%%%%%%%%%%%%%%%%%%%%%%%%%%%%%%%%%%%%%%%%%%%%%%%%%%%%%%%%%%%%%%%%%%%%%%%%%%%%%%%%%%%%%%%%%%%%%%%%%%%%%%%%%%%%%%%%%%%%%%%%%%%%%%%%%%%%%%%%%%%%%%%%%%%%%%%%%%%%%%%%%%%%%%%%%%%%%%%%%%%%%%%

\end{document}